\documentclass[twosided,12pt]{article}
\usepackage{amsfonts}

\usepackage{amscd}
\usepackage{amstex}

\newtheorem{theorem}{Theorem}[section]
\newtheorem{lemma}[theorem]{Lemma}
\newtheorem{proposition}[theorem]{Proposition}
\newtheorem{definition}[theorem]{Definition}
\newtheorem{remark}[theorem]{Remark}

\begin{document}

\author{J\o rgen Ellegaard Andersen and Josef Mattes}
\title{Configuration Space Integrals and Universal Vassiliev Invariants over Closed
Surfaces}
\date{28. April 1997}
\maketitle

\begin{abstract}
We show the existence of a universal Vassiliev invariant for links in
closed surface cylinders by explicit construction using configuration space
integrals. 
\end{abstract}

{\small \tableofcontents}

\pagestyle{myheadings}
\markboth{{\em J{\o}rgen Ellegaard Andersen and Josef
Mattes}}{{\em Configuration space integrals}}

\section{Introduction}

The quantization of the algebra of functions on the moduli space of flat
connections on a Riemann surface has attracted the attention of a number of
authors, including \cite{Alekseev95a}, \cite{Alekseev96b}, \cite{atiyah90}, 
\cite{Axelrod91}, [Buffenoir-Roche], \cite{Faltings93} and \cite{Hitchin90}.
The present authors together with Nicolai Reshetikhin developed a new
approach in \cite{Mattes96c}, \cite{Mattes96b},\cite{Mattes96a}. We showed
that chord diagrams (see \cite{Bar-Natan92}, \cite{Vassiliev90}) can be
generalized to surfaces and that they form a Poisson algebra of which the
algebra of functions on moduli space is a quotient in a natural way.
Furthermore, universal Vassiliev invariants (\cite{Kontsevich},\cite
{Bar-Natan96},\cite{Mattes96a}) allow us to use a simple, geometrically
defined multiplication of links to deformation quantize the algebra of chord
diagrams. In addition we showed that for punctured surfaces universal
Vassiliev invariants exist and that the quantization obtained descends to a
quantization of moduli space.

In the present paper we take this program one step further by constructing a
universal Vassiliev link invariant for closed surface cylinders. See definition
\ref{chord diagram} and theorem \ref{MT} for the statement of our main result. The
techniques are rather different from those we used before: For punctured
surfaces we used the combinatorial construction of universal invariants of
tangles due to Bar-Natan \cite{Bar-Natan93} and Le and Murakami \cite{Le95c}
whereas here we generalize the configuration space integral approach of Bott
and Taubes \cite{Bott}, D.Thurston \cite{Thurston95a}, Altschuler and\
Freidel \cite{Altschuler96a}, see also Guadagnini et al. \cite{Guadagnini90a}
and Bar-Natan \cite{Bar-Natan95a} as well as Axelrod and Singer \cite
{Axelrod94a}. Our constructions for punctured and closed surfaces
respectively are mutually orthogonal in a literal sense: In our construction
for punctured surfaces chords run parallel to the surface whereas in the
present construction chords run (almost) vertical.
In a subsequent paper we will show that the construction is local on the surface and study the
resulting deformation quantization of the algebra of functions on moduli space in detail.

We should point out that we are going to continue to use the notation
of \cite{Mattes96a} which differs from that preferred by some other
authors: We use $ch$ (rather than $\mathcal{A}$) to denote the space of
chord diagrams and $V$ (rather than $Z$) to denote a universal Vassiliev
link invariant. We also use the terms ``Vassiliev invariant'' and ``finite
type invariant'' interchangeably.

We wish to thank N.Reshetikhin and D.Thurston and to acknowledge that we
profited from a series of lectures by R.Forman at MSRI. Many of our
arguments generalize \cite{Thurston95a} and \cite{Altschuler96a} (which in turn are based on
\cite{Bott}).

\section{Diagrams in 3-manifolds}

We refer to \cite{Bar-Natan92} for basics concerning chord diagrams and
related matters. Let $M=M^{3}$ be an oriented 3-manifold.

\begin{definition}\label{chord diagram}
An \emph{abstract chord diagram} $D^{a}$ is a graph consisting of disjoint
oriented circles $L_{I}^{a},I\in \{1,...,n\}$ and disjoint arcs $
C_{j}^{a},j\in \{1,...,m\}$ such that:

\begin{enumerate}
\item  the endpoints of the arcs are distinct

\item  $\cup _{j}\partial C_{j}^{a}=\left( \cup _{i}L_{I}^{a}\right) \cap
\left( \cup _{j}C_{j}^{a}\right) $
\end{enumerate}

The arcs are called \emph{chords}, the circles are called the \emph{core
components} of the diagram. Set $L^{a}=\cup _{i}L_{I}^{a}$ and $C^{a}=\cup
_{j}C_{j}^{a}$.
\end{definition}

\begin{definition}
A \emph{chord diagram} $D$ in $M^{3}$ is a homotopy class of maps from an
abstract chord diagram $D^{a}$ into $M$. The space generated by chord
diagrams modulo 4T-relations is denoted by $ch\left( M\right) $.
\end{definition}

More generally we consider trivalent graphs $K^{a}$ containing a disjoint
collection of distinguished oriented circles $\left\{ L_{I}^{a}\right\} $.
Set $L^{a}=\cup _{i}L_{I}^{a}$ and $\Gamma ^{a}=\overline{K^{a}\backslash
L^{a}}$.

\begin{definition}
An \emph{internal vertex} of $K^{a}$ is a trivalent vertex of $\Gamma ^{a}$,
an \emph{external vertex} is a univalent vertex of $\Gamma ^{a}$. An \emph{
abstract Feynman diagram} is a trivalent graph $K^{a}$ together with a
cyclic orientation at each internal vertex such that every connected
component of $\Gamma ^{a}$ has at least one external vertex.
\end{definition}

If $e$ denotes the number of edges of $\Gamma ^{a}$, $u$ the number of
univalent vertices of $\Gamma ^{a}$, $u_{I}$ the number of univalent
vertices on $L_{I}^{a}$ (so that $u=\sum u_{I}$) and $t$ the number of
trivalent vertices, then $e=\frac{u+3t}{2}$. We define the degree of $K^{a}$
to be 
\[
\deg K^{a}=\frac{u+t}{2}=e-t
\]
Note that $t\leq 2\deg K^{a}$ and thus $e\leq 3\deg K^{a}$. For a chord
diagram the degree is just the number of chords. We also use the letter $e$ to
denote a general edge of $K^a$. It should be clear from the context, what the
meaning is in a given situation.

We denote the set of abstract Feynman diagrams with $u$ univalent and $t$
trivalent vertices by $\mathcal{F}_{u,t}^{a}$ and the set of abstract
Feynman diagrams of degree $n$ by $\mathcal{F}_{n}^{a}$. 

\begin{definition}
A {\em Feynman diagram}  is a homotopy class of mappings $K:K^{a}\rightarrow M^{3}$ 
from an abstract Feynman diagram $K^a$ into $M$, such that $\Gamma =K\left( \Gamma
^{a}\right) \subseteq M$ is contractible in $M$. We call $\Gamma =K\left(\Gamma ^{a}\right)
$ the \emph{graph} of the Feynman diagram and $
L_{K}=K\left( L^{a}\right) $ the \emph{core} of $K$.
\end{definition}

\begin{definition}
The vector space generated by Feynman diagrams in $M$ is denoted by $
\mathcal{F}\left( M\right) $.
\end{definition}

\begin{definition}
An \emph{automorphism} of a Feynman diagram $K$ is an automorphism of the
abstract Feynman diagram that can be realized by a homotopy of $K$ in $M$.
\end{definition}

The following generalizes a basic result for $\Bbb{R}^{3}$:

\begin{lemma}
We have that

\begin{enumerate}
\item  \label{a260197}$ch\left( M\right) $ $\cong $ $\mathcal{F}\left(
M\right) /STU$

\item  $\mathcal{F}\left( M\right) /STU$ satisfies AS and IHX.
\end{enumerate}
\end{lemma}

\noindent\textbf{Proof: }
The proofs of \cite[sect.3.1]{Bar-Natan92} apply since $\Gamma $ is
contractible by definition.
$\square$

\begin{remark}
Note that we have to be careful when applying STU: Two of the diagrams in
the STU relation having contractible graph does not imply that the third
diagram also has contractible graph.
\end{remark}

\section{Configuration spaces}
Throughtout this section $M$ will be an oriented compact $3$-manifold, possibly with
boundary. We let $M^{o}$ denote the interior of $M$.

\subsection{Definition and basic properties}

\label{b010407}

Let us briefly recall the basic constructions from \cite{Fulton94a}, \cite
{Bott}, \cite{Axelrod94a},\cite{Thurston95a}:

Let $\mathcal{C}_{r}(M)$ be the functorial compactified configuration space
of $r$ points in $M$ as described by Axelrod and Singer in \cite{Axelrod94a}:  Let
$[r]=\{1,\ldots ,r\}$. For any subset $S\subset [r]$ (with $|S|>1$,
which we will tacitly assume whenever necessary) we define the
diagonal in $M^{r}$ corresponding to $S$ by
\[
\Delta _{S}=\left\{ \left( x_{1},\dots ,x_{r}\right) |x_{I}=x_{j}\text{ for }
i,j\in S\right\} \subseteq M^{\times r}.
\]
Now we set 
\[
\mathcal{C}_{r}^{o}(M)=M^{r}\backslash \left( \cup _{S}\Delta _{S}\right) 
\]
and define the space $\mathcal{C}_{r}(M)$ to be the closure 
\[
\mathcal{C}_{r}(M)=\overline{\mathcal{C}_{r}^{o}(M)}\subset M^{r}\times
\left( \prod_{S\subset [r]}Bl_{S}\left( M\right) \right)
\]
where $Bl_{S}\left( M\right) $ is the (differential-geometric) blow-up of $
M^{r}$ along the diagonal $\Delta _{S}$, obtained by replacing the
diagonal $\Delta _{S}$ with the unit sphere bundle of the normal bundle to $
\Delta _{S}$. Equivalently, $Bl_{S}\left( M\right) =\left( M^{r}\backslash
\Delta _{S}\right) \cup \left( N\left( \Delta _{S}\right) \backslash\Bbb{R}^+\right)$ since 
\begin{equation}\label{a270497}
N\left( \Delta _{S}\right) \backslash\Bbb{R}^+ \cong \left\{ v\in TM^{\oplus S}:|v|^2 =1,
\sum_{I\in S} v_i =0 \right\}\times M^{\times ([r]-S) }
\end{equation}
By definition of $\mathcal{C}_r (M)$ we have projections from $\mathcal{C}_r (M)$ onto
$M^{\times r}$ and $Bl_{S}\left( M\right)$ for $S\subseteq [r]$. Hence, if $x\in \mathcal{C}_r
(M)$ is such that its image in $M^{\times r}$ is contained in $\Delta_S$ then the projection of
$x$ to $Bl_{S}\left( M\right)$ will be contained in $ N\left( \Delta _{S}\right)
\backslash\Bbb{R}^+$ and we simply write $(x_i )_{I\in S}$ for the element in $TM^{\oplus
S}$ which under the above isomorphism corresponds to the projection of $x$ to  $Bl_{S}\left(
M\right)$.

The space $\mathcal{C}_r (M)$ is a manifold with corners. If $M$ has no boundary, its 
strata $\mathcal{C}_r^{\mathcal{S}} (M)$ are parametrized by admissible collections of 
subsets of $[r]$: A collection $\mathcal{S}$ of subsets of $[r]$ is called {\em admissible} is and
$T\in \mathcal{S}$ has at least two elements and any two sets in $\mathcal{S}$
are either disjoint or one contains the other. The codimension of the
stratum is $\left| \mathcal{S}\right| $, in particular the codimension one
strata are parametrized by subsets $T$ of $\left[ r\right] $ of
cardinality at least two. 
If $M$ has a nonempty boundary, then there are more strata, since constraining
$p$ of the $r$ points to be contained in the boundary has codimension $p$.

The appropriate relative configuration space for a link $L$ with components $
L_{1},\dots ,L_{l}$ in $M$ is $\mathcal{C}_{t}\left( M,\left(
L_{1},u_{1}\right) ,\dots ,\left( L_{l},u_{l}\right) \right) $, introduced
in \cite[Appendix]{Bott} as the pullback 
\begin{equation}\label{a}
\begin{CD}
\mathcal{C}_{u,t}(L,M) @>>> \mathcal{C}_{u+t}(M) \\
@VVV                                            @VVV             \\
\times_i \mathcal{C}_{u_i}(L) @>>> \mathcal{C}_{u}(M)
\end{CD}
\end{equation}
where $\mathcal{C}_{r}\left( L\right) $ is defined just as $\mathcal{C}
_{r}\left( M\right) $ above and we used the following notation:
\[
\mathcal{C}_{u,t}\left( L,M\right) =\mathcal{C}_{t}\left( M,\left(
L_{1},u_{1}\right) ,\dots ,\left( L_{l},u_{l}\right) \right) 
\]
Let $u=\sum_{I=1}^{l}u_{I}$ where $\Gamma $ has $u_{I}$ univalent vertices
incident on $L_{I}$. Consider 
\[
\mathcal{C}_{t}\left( M,\left( L_{1},u_{1}\right) ,\dots ,\left(
L_{l},u_{l}\right) \right) ^{o}=\mathcal{C}_{t}\left( M,\left(
L_{1},u_{1}\right) ,\dots ,\left( L_{l},u_{l}\right) \right) \cap \mathcal{C}
_{t+u}^{0}(M)
\]

\begin{lemma}
The space $\mathcal{C}_{t}\left( M,\left( L_{1},u_{1}\right) ,\dots ,\left(
L_{l},u_{l}\right) \right) $ is the closure in $\mathcal{C}_{t+u}(M)$ of the
open configuration space $\mathcal{C}_{t}\left( M,\left( L_{1},u_{1}\right)
,\dots ,\left( L_{l},u_{l}\right) \right) ^{o}$ .
\end{lemma}

\noindent\textbf{Proof: }
The closure has the universal property of the pullback. $\square$

The stratification is given by 
\[
\mathcal{C}_{t}\left( M,\left( L_{1},u_{1}\right) ,\dots ,\left(
L_{l},u_{l}\right) \right) ^{\mathcal{S}}=\mathcal{C}_{t}\left( M,\left(
L_{1},u_{1}\right) ,\dots ,\left( L_{l},u_{l}\right) \right) \cap \mathcal{C}
_{t+u}^{\mathcal{S}} 
\]
for admissible $\mathcal{S}\subseteq \mathcal{P}\left( [u+t] \right) $. In
particular the codimension one boundary is given by 
\[
\partial ^{1}\mathcal{C}_{t}\left( M,\left( L_{1},u_{1}\right) ,\dots
,\left( L_{l},u_{l}\right) \right) =\bigcup _{T\subset [u+t]}\mathcal{C}
_{t}\left( M,\left( L_{1},u_{1}\right) ,\dots ,\left( L_{l},u_{l}\right)
\right) ^{\left\{ T\right\} } 
\]
and 
\begin{eqnarray*}
\partial ^{1}\left(\mathcal{C}_{t}\left( M,\left( L_{1},u_{1}\right) ,\dots
,\left( L_{l},u_{l}\right) \right) ^{\left\{ T\right\} }\right) &=&\bigcup
_{S\subset T}\mathcal{C}_{t}\left( M,\left( L_{1},u_{1}\right) ,\dots
,\left( L_{l},u_{l}\right) \right) ^{\left\{ T,S\right\} } \\
&&\bigcup _{S\supset T}\mathcal{C}_{t}\left( M,\left( L_{1},u_{1}\right)
,\dots ,\left( L_{l},u_{l}\right) \right) ^{\left\{ T,S\right\} } \\
&&\bigcup _{S\cap T=\emptyset }\mathcal{C}_{t}\left( M,\left(
L_{1},u_{1}\right) ,\dots ,\left( L_{l},u_{l}\right) \right) ^{\left\{
T,S\right\} }
\end{eqnarray*}

Let $\mathcal{C}_{u,t}(M)$ be defined as the space fibering over the space of links in $M$,
with fiber over $L$ equal $\mathcal{C}_{u,t}\left( L,M\right) $.

\begin{proposition}
The space $\mathcal{C}_{u,t}(M)$ is a manifold with corners. Its codimension one faces
are the codimension one faces of $\mathcal{C}_{u,t}(M^{o})$ together with $
\partial M \times \mathcal{C}
_{u,t-1}(M^{o})$.
\end{proposition}

\subsection{Correction bundle}

In analogy with \cite{Bott} and \cite{Thurston95a} we will introduce a
certain correction related to anomalous faces. To do this we construct for
each $T\subset [u+t]$ a bundle over the unit sphere bundle of $TM$ 
\[
\pi _{T}:\mathcal{B}_{T}\rightarrow S(TM).
\]
Let $T_{u}=T\cap [u]$. For a point $x$ in $\mathcal{C}_{u+t}^{\left\{
T\right\} }\left( M\right) $, we have $(x_i)_{I\in T}\in TM^{\oplus T}$ and we
define 
\[
D(x_i)=\text{span}\left\{ x_{I}\right\} 
\]
Likewise, for $v\in S(TM)$ we define $D(v)$ to be the line 
$\text{span}\left\{ v\right\} $ through $v$.

\begin{definition}
We set 
\[
\mathcal{B}_{T}^{o}=\left\{ (v,x)\in S(TM)\times \mathcal{C}_{u+t}^{\left\{
T\right\} }\left( M\right) |D(x_i)\subseteq D(v)\mbox{, } I\in T_u\right\} 
\]
and 
\[
\mathcal{B}_{T}=\overline{\mathcal{B}_{T}^{o}}\subseteq \mathcal{C}_{u+t}^{\left\{
T\right\} }\left( M\right) \subseteq \mathcal{C}_{u+t}\left(
M\right) 
\]
\end{definition}

For $v\in S(TM)$ we let $B_T(v)$ denote the fiber $\pi_T^{-1}(v)$.
The space $B_T$ is constructed such that we get the commutative diagram of subsets

\[
\begin{CD}
\mathcal{B}_{T} @>>> \mathcal{C}_{u+t}(M) \\
@AAA                                  @AAA                 \\
\mathcal{C}_{u,t}^{\{ T\} }(L,M) @>>> \mathcal{C}_{u,t}(L,M) \\
\end{CD}
\]

\begin{proposition}
The space $\mathcal{B}_{T}$ is a manifolds with corners of dimension 
\[
\dim (\mathcal{B}_{T})=3(u+t+1)-2|T_{u}|.
\]
The projection map 
\[
\pi _{T}:\mathcal{B}_{T}\rightarrow S(TM)
\]
is a fibration of manifolds with corners.
\end{proposition}

\begin{lemma}
\label{b110497}The codimension one boundary of $\mathcal{B}_{T}$ is given by 
$\partial ^{1}\mathcal{B}_{T}=\cup _{S}\mathcal{B}_{T}^{S}$ where $\mathcal{B
}_{T}^{S}=\mathcal{B}_{T}\cap \mathcal{C}_{t+u}^{\left\{ T,S\right\}(M) }$ and $
S\supset T,S\subset T$ or $S\cap T=\emptyset $.
\end{lemma}

Explicitly, $\mathcal{B}_{T}^{S}=\left\{ (v,x)\in S(TM)\times \mathcal{C}_{u+t}^{\left\{
S,T\right\} }|D(x_i)\subseteq D(v)\mbox{, } I \in T_u\right\} $

\subsection{Diagrams and strata}

\label{a240497}

Let $L$ be a link in $M$ and let $K^{a}\in \mathcal{F}_{u,t}^{a}$ be an
abstract Feynman diagram with $l$ core components, where $l$ is the number
of components of $L$. Choose an orientation of each of the edges in $\Gamma
^{a}$. Also choose a bijection $o$ between the set of components of $L$ and
the set of core components of $K^{a}$. Let $\{u_{1},\ldots ,u_{l}\}$,
respectively, be the number of univalent vertices on the core components of $
K^{a}$ so that $\sum_{I=1}^{l}u_{I}=u$.

Now choose one of each of the following four\textbf{\ }orderings:

\begin{itemize}
\item[$o_c$]  An ordering of the components of $L$.

\item[$o_u$]  An ordering of the $u_{I}$ univalent vertices on the $I$'th
core component of $K^{a}$ compatible with the cyclic order induced from the
orientation of the core component.

\item[$o_t$]  An ordering of the $t$ trivalent vertices of $K^{a}$.

\item[$o_e$]  An ordering of the $e$ edges of $\Gamma ^{a}$.
\end{itemize}

The orderings $o_{c}$ and $o_{u}$ together with the map $o$
specify a component of $\mathcal{C}_{u,t}(L,M)$. 

A collection $T$ of vertices of $\Gamma ^{a}$ defines a subgraph $\Gamma
^a_{T}\subseteq \Gamma ^{a}$ by $\Gamma ^a_{T}=\left\{ e|\partial e\subseteq
T\right\} $. Denote by $F^{d}$ the collection of those $
T\subseteq \left[ u+t\right] $ with $\Gamma ^a_{T}$ disconnected,  by $F^{a}$ the collection
of those $T\subseteq \left[
u+t\right] $ with $\left| T\right| >1$ and $\Gamma ^a_{T}$ is a connected
component of $\Gamma ^{a}$ and finally by $F^{h}$
the collection of those $T\subseteq \left[ u+t\right] $ with $\left|
T\right| >1$ and $\Gamma ^a_{T}$ connected but not a whole connected component
of $\Gamma ^{a}$ nor is $\Gamma^a_T$ just an edge of $\Gamma$, $F^p$ the
collection of those $T\subseteq \left[ u+t\right] $ for which there is an edge
in $K^a$ such that $\partial e = T$. 

We have the following decomposition of the codimension one strata of $
\mathcal{C}_{u,t}(M)$: 
\[
\partial \mathcal{C}_{u,t}^{1}(L,M)=\mathcal{C}_{u,t}^{b}(L,M)\cup
\mathcal{C}_{u,t}^{d}(L,M)\cup \mathcal{C}_{u,t}^{p}(L,M)\cup
\mathcal{C}_{u,t}^{h}(L,M)\cup 
\mathcal{C}_{u,t}^{a}(L,M),
\]
where 
\begin{eqnarray*}
\mathcal{C}_{u,t}^{b}(L,M) &=&\partial M
\times \mathcal{C}_{u,t-1}(L,M^{o}) \\
\mathcal{C}_{u,t}^{p}(L,M) &=&\cup _{T=F^p}\mathcal{C}_{u,t}^{\left\{
T\right\} }(L,M) \\
\mathcal{C}_{u,t}^{d}(L,M) &=&\cup _{T\in F^{d}}\mathcal{C}_{u,t}^{\left\{
T\right\} }(L,M) \\
\mathcal{C}_{u,t}^{h}(L,M) &=&\cup _{T\in F^{h}}\mathcal{C}_{u,t}^{\left\{
T\right\} }(L,M) \\
\mathcal{C}_{u,t}^{a}(L,M) &=&\cup _{T\in F^{a}}\mathcal{C}_{u,t}^{\left\{
T\right\} }(L,M)
\end{eqnarray*}
are the boundary, disconnected, principal, non-anomalous hidden and anomalous
strata, respectively (extending the terminology introduced by Bott and
Taubes in\textbf{\ }\cite{Bott}). By abuse of notation we will also call $T$
boundary, disconnected, principal, hidden or anomalous, respectively.

\subsection{Orientations}

\label{a010497}

We define the orientation of $\mathcal{C}_{u,t}(M)$ via an orientation of
the abstract Feynman graph, following \cite{Altschuler96a}. Let  $L$ be a link and let $K^{a}$
be an abstract Feynman diagram. 

\begin{proposition}
\label{a020297}There is an orientation on $\mathcal{C}_{u,t}$ depending on
the product of vertex orientations and product of edge orientations of $K^{a}
$.
\end{proposition}

\noindent\textbf{Proof: }
At a point $x\in \mathcal{C}_{u,t}^{o}(L,M)$, we
have that 
\[
T_{x}\mathcal{C}_{u,t}^{o}(L,M)\cong T_{x_1}L\oplus
\ldots \oplus T_{x_u}L\oplus T_{x_{u+1}}M\oplus
\ldots \oplus T_{x_{u+t}}M.
\]
Let $X_{x_{I}}=(X_{x_{I}}^{1},X_{x_{I}}^{2},X_{x_{I}}^{3})$ be an oriented
basis for $T_{x_{I}}M$, $I=u+1,\ldots ,u+t,$ and let $X_{x_{I}}^{1}$ be an
oriented basis for $T_{x_i}L$, $I=1,\ldots ,u$. The orientation element 
\[
\Omega (K^{a})=\bigwedge_{e\in \Gamma ^{a}}\Omega _{e},
\]
where 
\[
\Omega _{e}=X_{v}^{o_{v}(e)}\wedge X_{u}^{o_{u}(e)}\quad \partial e=(u,v),
\]
and $o_{v}$ is the cyclic order at the internal vertex $v$, defines an
orientation on $T_{x}\mathcal{C}_{u,t}^{o}(L,M)$, which orients
$\mathcal{C}_{u,t}^{o}(L,M)$. We then orient all the lower strata of
$\mathcal{C}_{u,t}(L,M)$ by inducing their orientation from $\mathcal{C}_{u,t}^{o}(L,M)$. 
$\square$

\subsection{Fibre integration}

\label{b260197}

We recall the following fact from \cite{Bott}: If $\pi :B\rightarrow X$ is a
fibration with fibre $F$ a compact manifold with corners then for any cycle $
c\subseteq X$ the value of the fibre integral $\pi _{*}(d\omega)
=\int_{F}d\omega $ of an exact form on $B$ is given by 
\begin{eqnarray*}
\left( \int_{F}d\omega \right) \left( c\right) &=&\int_{\pi ^{-1}c}d\omega
=\int_{\partial \pi ^{-1}c}\omega \\
&=&\int_{\pi ^{-1}\partial c}\omega +\int_{\pi ^{-1}c\cap \partial B}\omega
\\
&=&\int_{F}\omega \left( \partial c\right) +\int_{c}\left( \int_{\partial
F}\omega \right) \\
&=&d\left( \int_{F}\omega \right) \left( c\right) +\left( \int_{\partial
F}\omega \right) \left( c\right)
\end{eqnarray*}
so that $d\left( \int_{F}\omega \right) =\int_{F}d\omega -\int_{\partial
F}\omega $. In particular, $d\omega =0$ and $\int_{\partial F}\omega =0$
imply $\int_{F}\omega $ is constant.

\section{Configuration space integrals}

Let $\Sigma $ be a closed oriented surface of positive genus. For the remainder
of this paper $M=\Sigma\times I$ where $I=[0,1]$. Let $\pi _{\Sigma }$ denote the projection
from $M$ onto $\Sigma 
$ and $\pi _{I}$ the projection from $M$ onto $I$.

We choose a hyperbolic metric $\rho $ on $\Sigma $ if the genus is not 1. In
the case where the genus of $\Sigma $ is $1$ we choose a flat metric $\rho $. For any two points
$z_{1},z_{2}\in \Sigma $ we denote by $\rho
(z_{1},z_{2})$ the length of the shortest geodesic from $z_{1}$ to $z_{2}$.
For $x,y\in M$, we use the notation $\rho (x,y)$ for $\rho (\pi _{\Sigma
}(x),\pi _{\Sigma }(y))$. On $M$ we will consider the product metric.

\subsection{The two-form on $\mathcal{C}_{2}(M)$}

Fix a natural number $n$ (in the following $n$ will be the degree of a chord
diagram) and a small positive real number $\varepsilon $, smaller than the
injectivity radius of the metric $\rho$.

Consider 
\[
U_{n}=\left\{ (x,y)\in M\times M\mid \rho (\pi _{\Sigma }(x),\pi _{\Sigma
}(y))<(\varepsilon /3n)^{2}\right\} . 
\]
Let $N$ be an open neighbourhood of the diagonal in $M\times M$ which
satisfies that 
\[
N\cap U_{n}=\{(x,y)\in M\times M\mid |\pi _{I}(x)-\pi _{I}(y)|<\varepsilon
/3n\}\cap U_{n}. 
\]

Clearly we can choose $N$ such that $\mathcal{C}_{2}(M)$ is diffeomorphic to
the complement of $N$ in $M\times M$, say by a diffeomorphism $\iota :
\mathcal{C}_{2}(M)\rightarrow M\times M-N$. Note also that $U_{n}\cap
(M\times M-N)$ is disconnected.

Since $\varepsilon $ is smaller than the injectivity radius of $\Sigma $ we
have for any $(x,y)\in U_{n}$ a unique geodesic $\gamma =\gamma _{x,y}$ with 
$\gamma \left( 0\right) =\pi _{\Sigma }(x)$ and $\gamma \left( 1\right) =\pi
_{\Sigma }(y)$. Using this, we define a map $g:U_{n}\rightarrow T\Sigma $ by 
\[
g(x,y)=\gamma _{x,y}^{\prime }(0)\in T_{\pi _{\Sigma }(x)}\Sigma . 
\]

Let $\Lambda _{n}$ be a closed 2-form on $T\Sigma $ with support contained in 
$\{v\in T\Sigma \mid |v|<(\varepsilon /n)^{2}\}$ which represents the Thom
class of $T\Sigma $, see e.g. \cite[ch.I.6]{Bott1982}.

Define $\tilde{\omega}=\tilde{\omega}_{n}$, by 
\[
\tilde{\omega}_{n}=\left\{ 
\begin{array}{r@{\quad\mbox{ if }\quad}l}
g^{*}(\Lambda _{n}) & \pi _{I}(x)>\pi _{I}(y) \\ 
-g^{*}(\Lambda _{n}) & \pi _{I}(x)<\pi _{I}(y)
\end{array}
\right. 
\]

\begin{definition}\label{omegan.def}
Let $\omega =\omega _{n}$ be the closed 2-form on $\mathcal{C}_{2}(M)$ given
by 
\[
\omega =\iota ^{*}\tilde{\omega}_{n}
\]
\end{definition}

\begin{remark}\label{b270497}
By construction, $\omega $ has support in $\iota ^{-1}(U_{n})$.
\end{remark}

\subsection{The configuration space integral}
\label{tcsi}

Let $L$ be a link in $M$. An abstract Feynman diagram $K^{a}$ together with
the choices specified in section \ref{a240497} defines a map 
\[
\phi _{K^{a}}:\mathcal{C}_{u,t}\left( L,M\right) \rightarrow \mathcal{C}
_{2}\left( M\right) ^{\times e}
\]
defined on the component of $\mathcal{C}_{u,t}\left( L,M\right) $ determined by
these choices.

The form $\omega$ on $\mathcal{C}_{2}(M)$ can be pulled back to a two-form on 
$\mathcal{C}_{2}(M)^{\times e}$ via any of the projections $\pi :\mathcal{C}_{2}(M)^{\times
e}\rightarrow 
\mathcal{C}_{2}(M)$. The wedge product of these pullbacks defines a new form $\omega^e \in
\Omega^{2e}(\mathcal{C}_{2}(M)^{\times e})$.  

Consider the form $\phi _{K^{a}}^{*}\omega^{e} $ on the configuration space
$\mathcal{C}_{u,t}\left(
L,M\right)$. For every point $c\in \mathcal{C}_{u,t}\left( L,M\right) $ contained in the support
of this form there is
naturally induced a Feynman diagram $K_{c}^{a}$ in $M$ given by mapping the
edges of $K^{a}$ to the unique geodesics  between the corresponding pairs of
points in $M$. (The uniqueness follows from the fact that any pair of
endpoints of an edge is contained in $U_{n}$.)

Let $S_{l}$ denote the set of possible choices for the map $o$ (so $S_{l}$
is in bijective correspondence with the symmetric group on $l$ elements).

\begin{definition}
\label{c230297}Let $\omega _{\mathcal{C}_{u,t}}^{K^{a}}\in \Omega ^{2e}(
\mathcal{C}_{u,t}\left( L,M\right) ,\mathcal{F}(M))$ be given by 
\[
\omega _{\mathcal{C}_{u,t}}^{K^{a}}(c):=\sum_{o\in S_{l}}\frac{\phi
_{K^{a}}^{*}\omega^e (c)}{\left| \text{Aut}K_{c}^{a}\right| }K_{c}^{a}
\]
for $c\in \mathcal{C}_{u,t}\left( L,M\right) $. 
\end{definition}

The reason why we only sum over the choices of $o$ will be explained below.
Clearly, the class of $K_{c}^{a}$ is locally constant and therefore constant
on the connected components of $\text{supp}(\phi _{K^{a}}^{*}\omega^e
)\subseteq \mathcal{C}_{u,t}(L,M)$. Thus $\omega _{\mathcal{C}_{u,t}}^{K^{a}}
$ is integrable, and by compactness we see that $\int_{\mathcal{C}_{u,t}}\omega
_{\mathcal{C}_{u,t}(L,M)}^{K^{a}}$ is a finite sum of Feynman
diagrams in $M$, hence it specifies an element in $\mathcal{F}(M)$.

Since $\omega $ is a two form, this integral is independent of the chosen
total ordering $o_{e}$ of the edges of $\Gamma ^{a}$. From $\omega _{\left(
y,x\right) }=-\omega _{\left( x,y\right) }$ it follows that $\omega _{
\mathcal{C}_{u,t}}^{K^{a}}$ depends on the product of the orientations of
the edges. Thus proposition \ref{a020297} implies that $\int_{\mathcal{C}_{u,t}}\omega
_{\mathcal{C}_{u,t}}^{K^{a}}$ is independent of the choice of
edge orientations. Changing the ordering $o_{t}$ of the internal verticies
of $K^{a}$ or the ordering of the external verticies $o_{u}$ on any of the
components, results in an orientation preserving self-diffeomorphism of $
\mathcal{C}_{u,t}(L,M)$ taking the form for one orientation to the form for
the other, hence the integral is independent of these orderings. 

Hence we can now make the following definition:

\begin{definition}
We define for any abstract trivalent Feynman diagram $K^{a}$ 
\[
\tilde{V}\left( L,K^{a}\right) =\int_{\mathcal{C}_{u,t}\left( L,M\right)
}\omega _{\mathcal{C}_{u,t}}^{K^{a}}
\]
\end{definition}

We will consider the variation of this integral with respect to the link $L$. Using arguments
similar to those of \cite{Bott}, \cite{Thurston95a} and 
\cite{Altschuler96a} we calculate  $\delta \int_{\mathcal{C}_{u,t}\left(
L,M\right) }\omega _{\mathcal{C}_{u,t}\left( L\right) }^{K^{a}}$ using the
properties of the form $\omega _{\deg K}$ chosen above and the definition of
the automorphism group of a diagram.

\subsection{Variation of the function $\tilde{V}$}

\label{g010497}

Let us define 
\[
\tilde{V}_{n}(L)=\sum_{K^{a}\in \mathcal{F}_{n}^{a}}\tilde{V}(L,K^{a})
\]
and 
\[
\tilde{V}(L)=\sum_{n=0}^{\infty }\frac{\tilde{V}_{n}(L)}{2^{n}}
\]
thought of as a function on the space of links in $M$ with values in $\overline{ch\left( M\right)
}$. We want to compute the derivative of the
function $\tilde{V}$.

By the version of Stokes theorem presented in section \ref{b260197}, 
\[
\delta \tilde{V}=\delta ^{b}\tilde{V}+\delta ^{d}\tilde{V}+\delta ^{p}\tilde{V}+\delta
^{h}\tilde{V}+\delta ^{a}\tilde{V}
\]
where 
\begin{eqnarray*}
\delta ^{b}\tilde{V} &=&\sum_{n=0}^{\infty }\frac{1}{2^{n}}\sum_{K^{a}\in 
\mathcal{F}_{n}^{a}}\int_{\mathcal{C}_{u,t}^{b}\left( L\right) }\omega _{
\mathcal{C}_{u,t}\left( L\right) }^{K^{a}} \\
\delta ^{d}\tilde{V} &=&\sum_{n=0}^{\infty }\frac{1}{2^{n}}\sum_{K^{a}\in 
\mathcal{F}_{n}^{a}}\int_{\mathcal{C}_{u,t}^{d}\left( L\right) }\omega _{
\mathcal{C}_{u,t}\left( L\right) }^{K^{a}} \\
\delta ^{p}\tilde{V} &=&\sum_{n=0}^{\infty }\frac{1}{2^{n}}\sum_{K^{a}\in 
\mathcal{F}_{n}^{a}}\int_{\mathcal{C}_{u,t}^{p}\left( L\right) }\omega _{
\mathcal{C}_{u,t}\left( L\right) }^{K^{a}} \\
\delta ^{h}\tilde{V} &=&\sum_{n=0}^{\infty }\frac{1}{2^{n}}\sum_{K^{a}\in 
\mathcal{F}_{n}^{a}}\int_{\mathcal{C}_{u,t}^{h}\left( L\right) }\omega _{
\mathcal{C}_{u,t}\left( L\right) }^{K^{a}} \\
\delta ^{a}\tilde{V} &=&\sum_{n=0}^{\infty }\frac{1}{2^{n}}\sum_{K^{a}\in 
\mathcal{F}_{n}^{a}}\int_{\mathcal{C}_{u,t}^{a}\left( L\right) }\omega _{
\mathcal{C}_{u,t}\left( L\right) }^{K^{a}}.
\end{eqnarray*}

\begin{theorem}
The variation of $\tilde{V}$ equals the variation over the anomalous faces: 
\[
\delta \tilde{V}=\delta ^{a}\tilde{V}
\]
\end{theorem}

\noindent\textbf{Proof: }
By the above, we need to establish that $\delta ^{b}\tilde{V}=0$, $\delta
^{d}\tilde{V}=0$, $\delta ^{p}\tilde{V}=0$ and $\delta ^{h}\tilde{V}=0$.

The variation over the boundary faces $\delta ^{b}\tilde{V}$ is zero. This
follows since $\omega _{\left( x,y\right) }$ depends on $\pi _{\Sigma
}\left( x\right) ,\pi _{\Sigma }\left( y\right) $ only up to sign. Hence 
\[
\int_{x\in \Sigma \times \left\{ 1\right\} }\int_{\mathcal{C}_{u,t-1}\left(
\Sigma \times \left( 0,1\right) \right) }\phi _{K^{a}}^{*}\omega^e
=-\int_{x\in \Sigma \times \left\{ 0\right\} }\int_{\mathcal{C}_{u,t-1}\left( \Sigma \times \left(
0,1\right) \right) }\phi
_{K^{a}}^{*}\omega^e
\]
which shows that these contributions cancel.

If $T$ is disconnected $\phi _{K^{a}}$ factors through a codimension three
subspace so that $\phi _{K^{a}}^{*}\omega^e =0$.

The variation of $\tilde{V}$ coming from the hidden faces vanishes: By the
symmetry arguments in \cite{Thurston95a} there are either orientation
reversing automorphisms of the strata preserving $\omega ^{K}$ or
orientation preservering automorphisms of the strata mapping $\omega ^{K}$
to $-\omega ^{K}$ (this uses $\omega _{\left( y,x\right) }=-\omega _{\left(
x,y\right) }$).

Finally, the following section will establish that $\delta ^{p}\tilde{V}=0$. 
$\square$

The variation on anomalous faces $\delta ^{a}\tilde{V}$ might be nonzero.
Let $T\subset [u+t]$ define an anomalous face. Let $I : L \rightarrow
S(TM)$ be the natural tangential inclusion map. Consider the diagram 
\begin{equation}\label{e010497}
\begin{CD}
\label{a220197}
( \mathcal{C}_2 (M))^e   @= ( \mathcal{C}_2 (M))^e   @= ( \mathcal{C}_2 (M))^e  \\
@AA\phi_{K^{a}}A       @AA\phi_{K^{a}}A                @AA\psi_{K^{a}}A      \\
\mathcal{C}_{u,t} @<<<  \mathcal{C}_{u,t}^{\{ T\} } @>>>  B_T   \\
@VVV                  @VVV                                      @VV\pi_T V                \\
L^{\times u} \times M^{\times t}  @<<<   L   @>I >> S(TM) 
\end{CD}
\end{equation}
We obtain a form $\omega ^{K^a}_{B_{T}}$ on $B_T$ by replacing $\phi_{K^a}$ by
$\psi_{K^a}$ in definition \ref{c230297}.
If we let $\Omega_T^{K^a}$ be the 2-form on $S(TM)$ obtained by pushing
$\omega_{B_T}^{K^a}$ forward to $S(TM)$, i.e.
$$\Omega_T^{K^a}(v) = \int_{B_T(v)} \omega_{B_T}^{K^a},$$ 
then clearly
\begin{equation}
\delta ^{a}\int_{\mathcal{C}_{u,t}}\omega _{\mathcal{C}_{u,t}}^{K^{a}}=\int_{
\mathcal{C}_{u,t}^{\left\{ T\right\} }}\omega
_{B_{T}}^{K^{a}}|_{\mathcal{C}_{u,t}^{\{T\}}}=
\int_{L}i^{*}\Omega ^{K^a}_T  \label{b240497}
\end{equation}
In section \ref{d240497} we will adjust $\tilde{V}$ to obtain a function
which also has zero variation over the anomalous faces, hence is an
invariant of links.

\subsection{Variation on principal faces}

\label{c010497}

Recall that we are using $e$ both as a variable running through the edges of $K^{a}$
and as the number of edges in  $\Gamma^a$.

Collapsing an edge $e$ of an abstract Feynman diagram $K^{a}$ leads to a
graph $\delta _{e}K^{a}$ with one 4-valent vertex moreover, $e$ specifies a
principal face  $\mathcal{C}^{\partial e}_{u,t}\left( L,M\right) $. It is easy to
extend the discussion in section \ref{tcsi} to this case (for details see \cite[prop.3]
{Altschuler96a}). In particular, one $\delta _{e}K^{a}$ induces an orientation
on $\mathcal{C}^{\partial e}_{u,t}\left( L,M\right) $ and one gets a map  $\phi _{\delta
_{e}K^{a}}
: \mathcal{C}^{\partial e}_{u,t}\left( L,M\right) \rightarrow C_2(M)^{\times e}$
such that  
\[
\omega^{K^a}_{{\mathcal C_{u,t}}}|_{\mathcal{C}_{u,t}^{e}(L,M)}(c)=
\sum_{o\in S_l}\frac{\phi _{\delta _{e}K^{a}}^{*}\omega^e
(c)}{\left| \text{Aut}K^a_{c}\right| }K_{c}^{a}
\]

\begin{proposition}
\label{b190297} The orientation defined in proposition \ref
{a020297} and the orientation induced from $\delta_e K^a$ on $\mathcal{C}^{\partial
e}_{u,t}\left( L,M\right) $
agree.
\end{proposition}

\noindent\textbf{Proof: }
As in \cite[prop.3]{Altschuler96a} 
$\square$

\begin{lemma}
Let $\Gamma _{x}$ be a graph with one four-valent vertex $x$, all other
vertices trivalent. Then there are at most three trivalent graphs $\Gamma
_{x}^{I},I\in \left\{ 1,2,3\right\} $, such that $\delta _{e_{I}}\Gamma
_{x}^{I}=\Gamma _{x}$ for edges $e_{I}\in \Gamma _{x}^{I}$.
\end{lemma}

\noindent\textbf{Proof: }
Holds as in \cite[lemma 1]{Altschuler96a} since this is a local statement.
$\square$

\begin{lemma}
$K$ a trivalent graph with edges $e,f$. Then $\delta _{e}K\cong \delta _{f}K$
iff there is an $\sigma \in \text{Aut}\left( K\right) $ such that $\delta
_{e}\left( \sigma \left( K\right) \right) =\delta _{e}\left( K\right)
,\sigma \left( e\right) =f$
\end{lemma}

\noindent\textbf{Proof: }
Obvious. (as \cite[lemma 2]{Altschuler96a})
$\square$

\begin{proposition}
$\delta _{e}K\cong \delta _{f}K$ iff $\sigma \left( K\right) =\pm \left(
\delta _{e}K\right) ^{I},\sigma \left( f\right) =e$ for some $I,\sigma \in 
\text{Aut}\left( K\right) $
\end{proposition}

\noindent\textbf{Proof: }
As \cite[prop.1]{Altschuler96a}.
$\square$

\begin{lemma}
If $e$ is admissible and collapsing $e$ does not result in two vertices
being connected by two different edges, then $P_{I}^{\pm }\left( K,e\right)
\cap P_{j}^{\pm }\left( K,e\right) =\emptyset $
\end{lemma}

\noindent\textbf{Proof: }
\cite[lemma 5]{Altschuler96a} is stronger.
$\square$

\begin{lemma}
$\#\left\{ f|\delta _{f}K\cong \pm \delta _{e}K\right\} =\#\left\{ I\right\}
\times \frac{\left| \text{Aut}_{+}K\right| }{\left| \text{Aut}
_{+}\delta _{e}K\right| }$
\end{lemma}

\noindent\textbf{Proof: }
Follows from the above lemmata, as in \cite[lemma 4]{Altschuler96a}. 
$\square$

{}From the above it now follows that

\begin{lemma}
The variation on principal faces is given by 
\[
\delta ^{p}\tilde{V}=\sum_{n=0}^{\infty }\frac{1}{2^{n}}\sum_{\stackrel{K^a\in
{\mathcal F}_n}{e\in K^a}}\frac{1}{2^{n}}\int_{\mathcal{C}^{\partial e}_{u,t}\left(
L,M\right) }\omega _{\mathcal{C}_{u,t}}^{\delta _{e}K^{a}}
\]
\end{lemma}

\begin{lemma}
If there is an orientation reversing automorphism of $K$ or $\delta _{e}K$
then $\left[ K\right] =0$ or $\int {}_{\mathcal{C}^{\partial e}_{u,t}\left( L,M\right)
}\omega _{\mathcal{C}_{u,t}}^{\delta _{e}K^{a}}=0$.
\end{lemma}

\noindent\textbf{Proof: }
Obvious.
$\square$

\begin{lemma}
One gets that $\int_{\mathcal{C}^{\partial e}_{u,t}\left( L,M\right) }\omega
_{\mathcal{C}_{u,t}}^{\delta _{e}K^{a}}=0$ if collapsing $e$ does result in two vertices
being connected by two different edges.
\end{lemma}

\noindent\textbf{Proof: }
$\phi _{K^{a}}$ factors through $\left( \mathcal{C}_{2}\left( M\right)
\right) ^{e-1}$, hence $\phi _{K^{a}}^{*}\omega^e =0$.
$\square$

\begin{theorem}
The variation on principal faces is zero: $\delta ^{p}\tilde{V}=0$.
\end{theorem}

\noindent\textbf{Proof: }
Follows from the last three lemmata as in \cite[theorem 1]{Altschuler96a}
since all the arguments are local. 
$\square$

\section{The invariant}

\subsection{Correction terms}

\label{b020297}
 
Here we will construct the correction factor to cancel the variation over
the anomalous faces. From equation (\ref{b240497}) we see that the we need a
factor whose variation is $\sum_{K^a\in{\mathcal F}^a_n}\sum_{T\in F^a(K^a)}
\int_{L}i^*\Omega_T^{K^a}$. To construct such a
correction we want to show that the 2-form $\Omega_n = \sum_{K^a\in{\mathcal
F}^a_n}\sum_{T\in
F^a(K^a)}\Omega_T^{K^a}$ is exact. We
will do this only on an open dense submanifold $S'(TM)$ of $S\left( TM\right) $ containing a
neighbourhood of $L$, but this suffices for our purposes.

Our first step is

\begin{proposition}
The 2-form $\Omega_n$ is closed.
\end{proposition}

\noindent\textbf{Proof: }
By lemma \ref{b110497},
$$\delta \Omega_n=\sum_{K^a\in{\mathcal F}^a_n}\sum_{T\in
F^a(K^a)}\sum_{S}\int_{B_T^S(v)} \omega
^{K^a}_{B_{T}}$$

where one of the following holds:

\begin{enumerate}
\item  $S$ defines a boundary face 

\item  $S\subset T$ 

\item  $S\supset T$ 

\item  $S\cap T=\emptyset $
\end{enumerate}

If $S$ is boundary the contribution is zero as above.

If $S\subset T$ then $S$ is principal or hidden and the contribution is zero.

If $S\supset T$ then $S$ is disconnected and does not contribute.

Hence we see that the contributions to $\delta \Omega_n$ all vanishes unless $S\cap T =
\emptyset$ and $S$ is 
anomalous.

Now we note that $S$ contains at least one univalent vertex $v_{0}$ (since $S$ is a connected
component of $\Gamma $) that in $
\mathcal{C}_{u,t}^{\left\{ T,S\right\} }\left( L,M\right) $ runs over the
whole tangent space of $T_{x}M$ if the points parametrized by $S$ collapse
at $x$. Let $v_1$ be the vertex which is connected to $v_0$ by a edge in
$\Gamma$. Reflection of $v_0$ in $v_1$ 
is an automorphism $\theta$ of $B_T^S$ such that $\theta^* \omega = -\omega $ and $\theta^*
\int_{B_T^S(v)} 
\omega = - \int_{B_T^S(v)} \omega$ implies $\int_{B_T^S(v)} \omega=0$.
This completes the proof that $\delta \Omega_n=0$.
$\square$

Now we proceed as follows: 

Let $\frak{n}$ be a section of $S(TM)$. We shall now 
restrict our self to links
$L\subset M$ such that $\frak{n}(M)\cap TL = \emptyset$ inside $S(TM)$. In other words
$\frak{n}$ 
induces a framing of $L$.

Choose a point $x\in \Sigma$ not on the projection of the link $L$ and let
$$S'(TM) = S(TM)-\left\{ \frak{n}(M)\cup S(TM)|_{\pi_\Sigma^{-1}(x)}\right\}.$$
By construction $S'(TM)$ is homotopy equivalent to $\Sigma-\{x\}$, hence
$H^2(S'(TM))= 0$. Hence, if we let $\Omega'_n = \Omega_n|_{S'(TM)},$ then
$\Omega'_n$ is exact. Now pick cycles $b_{I}$ representing a basis for $
H_{1}\Sigma $.

\begin{definition}
\label{c110497}Let $\alpha _n$ be such that 
\begin{eqnarray*}
\Omega'_n &=&d\alpha _n \\
\alpha _n\left( b_{I}\right)  &=&0.
\end{eqnarray*}
\end{definition}

Any two choices for $\alpha _n$ will be cohomologous: $d\left( \alpha
_{n}-\alpha _{n}^{\prime }\right) =\Omega_n -\Omega_n =0$ and $\left( \alpha
_{n}-\alpha _{n}^{\prime }\right) \left( b_{I}\right) =0\forall I$ so that $
\left[ \alpha _{n}-\alpha _{n}^{\prime }\right] =0\in H^{1}\left( S'(TM)\right) $.

By the version of Stokes theorem in section \ref{b260197} we see that 
\begin{equation}
\delta \int_{L}\alpha _{n}=\int_{L}d\alpha _{n}=\int_{L}\Omega_n  \label{a150497}
\end{equation}
and $\int_{L}\left( \alpha _n-\alpha _{n}^{\prime }\right)
=\int_{L}df_{n}=\int_{\partial L}f_{n}=0$ shows that it is
independent of $\alpha _{T}$. Since $\frak{n}$ is only allowed to vary
continuously in the complement of $I(L)$ and since $
H^{2}\left( \Sigma \backslash \left\{ x\right\} \right) =0$, it is clear that the
cohomology class of $\alpha_n $ is independent of the choice of framing $\frak{
n}$.

Furthermore $\int_{L}\alpha_n $ is independent of the choice of $x$: Consider
two choices of points $x,y$ and the corresponding forms $\alpha ^x_{n},\alpha
^y_{n}$. The first homology of $S\left( TM\right) \backslash \left\{ \frak{n}(M)\cup 
S(TM)|_{\pi_\Sigma^{-1}(\{x,y\})}\right\}$ is generated by 
$\left\{ b_{I}\right\} $ and a small loop $c_{y}$ around $y$ on $\Sigma$. 
By assumption $\forall \alpha ^x_{n}\left( b_{I}\right)
=\alpha ^y_{n}\left( b_{I}\right) $ for all $I$, in addition $\alpha
^x_{n}\left( c_{y}\right) =0$ since $c_y$ is homologous to $0$ in the domain of $
\alpha ^x_{n}$. Since also $\alpha ^y_{n}\left( c_{y}\right) =\alpha ^y_{n}\left(
c_{x}\right) =0$ (because $c_{y}$ is homologous to $c_{x}\text{ mod}
\left\{ b_{I}\right\} $) we see that $\alpha ^x_{n}$ is cohomologous to $
\alpha ^y_{n}$ near $L$ and thus $\int_{L}\alpha ^x_{n}=\int_{L}\alpha _n^{y}$.

\subsection{Definition of the invariant}

\label{d240497}

Let $\overline{ch\left( M\right) }$ denote the completion of $ch\left(
M\right) $ with respect to the filtration defined by degree.

\begin{definition}
Let $V_n\left( L\right) $ be defined by 
$$V_n\left( L\right) =\tilde{V}_n(L) - \int_{L}\alpha_n.$$
\end{definition}

\begin{theorem}\label{MT}
The sum 
\[
V\left( L\right) =\sum_{n}\frac{1}{2^{n}}V_n(L)
\]
defines an element of $\overline{ch\left( M\right) }$. The map $L\mapsto
V\left( L\right) \in \overline{ch\left( M\right) }$ is a universal finite
type invariant for links in $M$.
\end{theorem}

\noindent\textbf{Proof: }

By the definition of $U_{n}$ we know that only those imbeddings of $K$
contribute where for any connected component $\bar{K}$ of $K$ the projection 
$\pi _{\Sigma }\left( \bar{K}\right) $ is contained in some disk of radius
at most $e\left( \frac{\varepsilon }{3\deg K}\right) ^{2}\leq \frac{
\varepsilon ^{2}}{3\deg K}$ on $\Sigma $, in particular the graph of $\Gamma 
$ is contractible. From this, the finiteness of the sum defining $V$ and from $\mathcal{F}
/STU\cong ch$ (lemma \ref{a260197}) it follows that the sum is well-defined
as a sum of chord diagrams.

We have that
$$\delta V_n = \delta \tilde{V}_n - \delta \int_{L} \alpha_n = \delta^a
\tilde{V}_n - \int_{L}\Omega_n = 0$$
by the construction of $\Omega_n$ in the previous section.

We are left with showing that $V\left( L\right) $ is in fact a universal
finite type invariant. This will be done in the next section.
$\square$

\subsection{Finite type and universality}

Consider a resolution of a singular knot $\sum_{\eta \in \left\{ \pm
1\right\} ^{n}}\left( -1\right) ^{\eta }L_{\eta }$ of degree $n$. Let the $
L_{\eta }$ be equal outside balls $B_{I},I=1,\dots ,n$, of size $\varepsilon 
$ around the singularities. By invariance of $V$ under isotopy we can also
assume they are almost planar (contained in $\Sigma \times \left\{ \frac{1}{2
}\right\} $), all the $\pi _{\Sigma }\left( L_{\eta }\right) $ are equal and
the distance of the balls satisfies $\left| B_{I},B_{j}\right| _{\mu
}>3\varepsilon /n$. If for some $I$ no edge of $\Gamma $ ends in $B_{I}$
then clearly $\sum \left( -1\right) ^{\eta }\omega _{\eta }^{K^{a}}=0$. 

Now we consider the correction term: Assume $L^{1}$ and $L^{2}$ are singular
links of degree $n-1$ that differ by a crossing switch away from the $B_{I}$
and let $\chi \left( t\right) $ be a homotopy from $L^{1}$ to $L^{2}$, then 
\begin{eqnarray*}
\sum_{\eta \in \left\{ \pm 1\right\} ^{n-1}}\left( -1\right) ^{\eta
}\int_{L_{\eta }^{1}}\alpha _{\deg K^{a}}-\sum_{\eta \in \left\{ \pm
1\right\} ^{n-1}}\left( -1\right) ^{\eta }\int_{L_{\eta }^{2}}\alpha _{\deg
K^{a}}\\
=\sum_{\eta \in \left\{ \pm 1\right\} ^{n-1}}\left( -1\right) ^{\eta
}\int_{t}\int \chi \left( t\right) ^{*}\Omega _{\deg K^{a}}^{\prime }.
\end{eqnarray*}
If $2\deg K^{a}<n-1$ then for each imbedding of $K^{a}$ either the
projection of some edge has length greater than $\left( \varepsilon
/3n\right) ^{2}$ or at least one $B_{I},I=1,\dots ,n-1$  does not contain a
univalent vertex. In either case $$\sum_{\eta \in \left\{ \pm 1\right\}
^{n-1}}\left( -1\right) ^{\eta }\int_{\text{Im}\chi \left( t\right) }\Omega
_{\deg K^{a}}^{\prime }=0$$ for each $t$, so that 
\begin{eqnarray*}
\lefteqn{\sum_{\eta ^{\prime }\in \left\{ \pm 1\right\} ^{n}}\left( -1\right) ^{\eta
^{\prime }}\int_{L_{\eta ^{\prime }}}\alpha _{\deg K^{a}}=} \\
& &\sum_{\eta \in
\left\{ \pm 1\right\} ^{n-1}}\left( -1\right) ^{\eta }\int_{L_{\eta
}^{1}}\alpha _{\deg K^{a}}-\sum_{\eta \in \left\{ \pm 1\right\}
^{n-1}}\left( -1\right) ^{\eta }\int_{L_{\eta }^{2}}\alpha _{\deg K^{a}}=0.
\end{eqnarray*}
This shows that $V\left( L\right) $ of finite type.

Finally, we show that $V\left( L\right) $ is universal: We observe that any
chord that contributes to $\int \omega ^{K^{a}}$ projects into a small disk
in $\Sigma $ (by remark \ref{b270497}) so that the integral decomposes as a
product of integrals over disjoint domains for each double point of the
projection since $\left| B_{I},B_{j}\right| _{\mu }>3\varepsilon /n\geq 
\frac{e\varepsilon }{n\deg K}$ (cf. \cite[proof of theorem 5]{Altschuler96a}
), hence if a component of $\Gamma $ intersects two different $B_{I}$ then
the contribution to the invariant vanishes. Let $\gamma $ be the number of
connected components of $\Gamma ^{a}$.  Since  $\gamma \leq \deg K^{a}$ and $
\gamma =\deg K^{a}$ iff $K^{a}$ has no trivalent vertex it follows that $
\sum_{\eta \in \left\{ \pm 1\right\} ^{n}}\left( -1\right) ^{\eta }\tilde{V}
_{I}\left( L_{\eta }\right) =0$ for $i<n$. Furthermore, $\alpha $ can
contribute only if $K^{a}$ has at least one trivalent vertex and $\sum_{\eta
\in \left\{ \pm 1\right\} ^{n-1}}\left( -1\right) ^{\eta }\int_{\text{Im}
\chi \left( t\right) }\Omega _{I}^{\prime }=0$ for $I=1,\dots ,n$ by the
argument above so that $\sum_{\eta \in \left\{ \pm 1\right\} ^{n}}\left(
-1\right) ^{\eta }\int_{L_{\eta }}\alpha _{I}=0$ for $i<n$.

If $e=\deg K^{a}=n$ (i.e. no trivalent vertex) then $\tilde{V}\left( L_{\eta
},K^{a}\right) \neq 0$ only if there is exactly one chord in each $B_{I}$,
which necessarily runs almost vertical. The coefficient will be $1$ by the
defining properties of the Thom class $\Lambda _{n}$ entering in the
construction of $\omega _{n}$. To see that the correction term will be zero
for $\deg K^{a}=n$ we note that we can assume that $\chi \left( t\right) $
moves only in the vertical direction. By the definition of $\omega $ again,
we see that $\Omega ^{\prime }$ has no component in the vertical direction,
thus we effectively integrate a two-form over a manifold of dimension one.
Hence the anomalous contribution is zero for $\deg K^{a}=n$.

This concludes the proof of the main theorem \ref{MT}.

\bibliographystyle{alpha}
\bibliography{books,endnote}

Our \textbf{preprints} are available from \emph{http://math.ucdavis.edu/
\symbol{126}mattes}

\textbf{Addresses:}

J.E.A.: Dept. of Math., Univ. of Aarhus, DK-8000 Aarhus C, DK

J.M.: Dept. of Math., UC Davis, CA-95616, USA

\end{document}